\documentclass{amsart}
\usepackage{graphicx,mathrsfs,amsmath,amssymb}
\newtheorem{Theorem}{Theorem}[section]

\newtheorem{Lemma}[Theorem]{Lemma}

\theoremstyle{Definition}
\newtheorem{Definition}[Theorem]{Definition}
\theoremstyle{Remark}
\newtheorem{Remark}[Theorem]{Remark}
\numberwithin{equation}{section}

\begin{document}

\title[Recovery of inhomogeneities and inclusions]{Recovery of inhomogeneities and buried obstacles}
\author{Hongyu Liu}
\address{University of Washington, Department of Mathematics, Box 354350, Seattle, WA 98195, USA}
\email{hyliu@math.washington.edu}

\thanks{The work is
partly supported by NSF grant, FRG  DMS 0554571.}

\begin{abstract}
In this paper we consider the unique determination of
inhomogeneities together with possible buried obstacles by
scattering measurements. Under the assumption that the buried
obstacles have only planar contacts with the inhomogeneities, we
prove that one can recover both of them by knowing the associated
scattering amplitude at a fixed energy.
\end{abstract}
\maketitle
\section{Introduction}

We shall be concerned with the unique determination of a medium
together the possible buried obstacles by making scattering
measurement far away from the (unknown/inaccessible) object. There
are no global identifiability results available in literature on
recovering both of them by knowing the scattering amplitude at a
fixed energy. The existing results are either based on knowing the
outside inhomogeneity in advance to recover only the included
obstacle (\cite{KirPai}); or by making use of measurement data
with frequency from an open interval (\cite{Hah}), that is, much
more data are utilized than needed. Moreover, it is noted that the
uniqueness result in \cite{Hah} cannot be generalized to three
dimensions since its argument involves conformally mapping the
domain containing the support of the inhomogeneity onto an
annulus. We would like to mention that the global uniqueness in
the determination of scatterers consisting of sole mediums or
obstacles by scattering amplitude at a fixed energy has been
widely known as sophisticatedly established. The result for
recovering a medium was obtained by Nachman (\cite{Na}) which is
based on the use of complex geometric optics (CGO) solutions due
to Sylvester-Uhlmann (\cite{SU}); and for recovering an obstacle
was obtained by Kirsch-Kress (\cite{KirKre}) which is based on the
use of singular sources due to Isakov (\cite{Isa}). We also refer
to the monographs \cite{ColKre} and \cite{Isa} for a comprehensive
discussion and related literature.

In this paper, in combination of the two methodologies of using
CGO solutions and singular sources we are able to prove the global
identifiability of both the scattering medium and the possible
buried obstacles. Our restrictive assumption is that the buried
obstacles have only planar contacts with the inhomogeneities. In
general, we cannot recover an obstacle which is completely buried
inside the inhomogeneity. But the exposure part of the obstacle to
the exterior of the medium can be arbitrarily small. On the other
hand, our proofs indicate that if the obstacle is enclosed
entirely in the medium but known in advance, then one can recover
the surrounding medium by the corresponding scattering amplitude.
In the rest of this section, we give a brief formulation of the
direct and inverse scattering problems.

Let $\mathbf{D}\subset \mathbb{R}^3$ represent the obstacle which
is a bounded domain with connected Lipschitz complement
$\mathbf{G}:= \mathbb{R}^3\backslash \bar{\mathbf{D}}$; that is,
we include in our discussion the case of multiple obstacle
components. Further let $\mathrm{B}$ be a sufficiently large ball
such that $\bar{\mathbf{D}}\subset \mathrm{B}$. Let $q\in
L^\infty(\mathbf{G})$ with $\Im q\geq 0$ and $supp (1-q)\subset
\mathrm{B}\backslash {\mathbf{D}}$ represent the scattering
medium, namely, the \emph{refractive index}. We consider the
following scattering problem for the time-harmonic plane wave
$u^i(x):=\exp\{i\kappa x\cdot \theta\}$
\begin{equation}\label{eq:direct system}
\begin{cases}
 \,\,(\Delta+q\kappa^2) u=0\, &\mbox{in $\mathbf{G}$},\\
 \quad \mathcal{B}u=0 &\mbox{on $\partial \mathbf{D}$},\\
 \quad \mathcal{M}u=0 &\mbox{in $\mathbb{R}^3$},
\end{cases}
\end{equation}
where $u:=u^i+u^s$ with $u^s$ the so-called scattered field and $\mathcal{B}$ is the boundary operator which
gives a Dirichlet boundary condition $u|_{\partial \mathbf{D}}=0$ corresponding to a \emph{sound-soft} obstacle $\mathbf{D}$.
Moreover, the last equation is the well-known \emph{Sommerfeld
radiation condition} given by
\begin{equation}\label{eq:sommerfeld}
\lim\limits_{r\rightarrow\infty}r(\frac{\partial u^s}{\partial
r}-i\kappa u^s)=0\quad r=|x|,
\end{equation}
which holds uniformly for all directions $\hat{x}:=x/|x|\in\mathbb{S}^2$ . It is known that $u$
has the following asymptotic representation (see \cite{ColKre})
\begin{equation}\label{eq:far field}
u(x;\theta,\kappa)=\exp\{i\kappa\theta\cdot x\}+\frac{\exp\{i\kappa|x|\}}{|x|}\mathcal{A}(\hat{x},\theta,\kappa)
+\mathcal{O}(|x|^{-2}).
\end{equation}
The function $\mathcal{A}$ is called the \emph{scattering amplitude} (or the \emph{far-field pattern}) with $\hat{x}$, $\theta$
and $\kappa$ denoting, respectively, the observation direction, the incident direction and the wave number.
The inverse scattering problem consists in the determination of the obstacle $\mathbf{D}$ and the scattering medium $q$ by knowing $\mathcal{A}(\hat{x},\theta,\kappa)$ for a fixed $\kappa>0$ and all $\hat{x}\in \mathbb{S}^2$, $\theta\in \mathbb{S}^2$.

The paper is organized as follows. In section~\ref{sect:direct
scattering}, we present the class of admissible scatterers and
give a brief study of the forward scattering problem.
Section~\ref{sect:uniq with soft obstacle} is devoted to the
unique determination of a scatterer with the buried obstacle. In
Section~\ref{sect:concluding remark}, we indicate how to determine
the surrounding medium when the obstacle is buried inside
completely but known \emph{a priori}.

\section{Class of admissible scatterers and the direct scattering
problem}\label{sect:direct scattering}

In order to state our uniqueness results we need first to introduce a suitable class $\mathcal{C}$ of admissible scatterers.
We begin by fixing some notations which shall be used throughout the rest of the paper. For any
$x\in \mathbb{R}^3$ and $r>0$, with $B_r(x)$ we denote the open ball of center $x$ and radius $r$. Let
$\Gamma_l$ with an index $l\in \mathbb{N}$ represent a simply connected subset of some plane $\Pi_l$
in $\mathbb{R}^3$, and moreover, $\mathscr{R}_l$ represent the reflection in $\mathbb{R}^3$ with respect to
$\Pi_l$. Let $C$ denote a generic constant which may be changed in different inequalities but must be fixed and finite
in a given relation. Finally, $``a\lesssim b"$ shall refer to $``a\leq C b"$.

\begin{Definition}\label{def:scatterer}
We say that $\Sigma(\mathbf{D},q,\Omega)$ is a scatterer of class
$\mathcal{C}$ with obstacle $\mathbf{D}$ and scattering medium $q$
if it satisfies the following assumptions
\medskip
\begin{enumerate}
\item[i)] $\Sigma$ is a compact set in $\mathbb{R}^3$ with connected complement $\Sigma_\infty:=\mathbb{R}^3
\backslash {\Sigma}$.
\medskip
\item[ii)] $\Sigma=\overline{\Omega\cup \mathbf{D}}$, where $\mathbf{D}$ is a $C^{2,1}$ domain with connected complement
$\mathbf{G}:=\mathbb{R}^3\backslash\bar{\mathbf{D}}$ and
$\Omega\subset \mathbf{G}$ is an open set.
\medskip
\item[iii)] $q(x)\in L^\infty(\mathbf{G})$ with $\Im q\geq 0$ and $supp(1-q)=\bar{\Omega}$. Moreover, $1-q\in C^{0,1}(\bar{\Omega})$
and there exists a constant $\epsilon_0>0$ such that $|1-q(x)|\geq \epsilon_0$ for all $x\in \bar{\Omega}$, i.e., $1-q$ has jump across
$\partial \Omega$.
\medskip
\item[iv)] The medium and the obstacle have only planar contact in the sense that $\partial \Omega\cap \partial \mathbf{D}=
\bigcup_{l=1}^{N_0}\Gamma_l$, where $N_0$ is a finite integer and $\Gamma_l\subset\partial \Omega_l$ with each
$\Omega_l$ a connected component of $\Omega$, $l=1,2,\ldots, N_0$. Moreover, if we set $\Omega_0:=\Omega-\bigcup_{l=1}^{N_0}\Omega_l$,
then $\bar{\Omega}_l\cap
\bar{\Omega}_{l'}=\emptyset$ if $l\neq l'$ for $0\leq l,l'\leq N_0$.
\medskip
\item[v)] Set $\Gamma_{int}:=\bigcup_{l=1}^{N_0}\Gamma_l$ and $\partial\mathbf{D}_{ext}:=\partial \mathbf{D}\backslash \Gamma_{int}$,
$\partial\Omega_{ext}:=\partial \Omega\backslash \Gamma_{int}$. $\partial\Omega_{ext}$ are $C^{2}$ continuous.
\medskip
\item[vi)] $(\Omega_l\cup
\mathscr{R}_l\Omega_l)\cap (\Omega_{l'}\cup
\mathscr{R}_{l'}\Omega_{l'})=\emptyset$ for $1\leq l, l'\leq N_0$
and $l\neq l'$; and $(\Omega_l\cup
\mathscr{R}_l\Omega_l)\cap \Omega_0=\emptyset$ for $1\leq l\leq N_0$.
\end{enumerate}
\end{Definition}


Clearly, according to our definition, a scatterer $\Sigma(\mathbf{D},q,\Omega)\in \mathcal{C}$ is composed of an
impenetrable obstacle $\mathbf{D}$ and the surrounding medium
$q$ with support in $\Omega$; and the $\Gamma_{int}$ part of the obstacle $\mathbf{D}$ is buried in the inhomogeneity.
Moreover, it is noted that we admit multiple scattering components. In fact, if we let $\sigma$ denote a connected component
of $\Sigma$, then it may be an obstacle, or the support of a scattering medium, or the two combined together with the
obstacle buried inside the inhomogeneity. Hence, an admissible scatterer is much general which may consist of
multiple components being obstacles, or scattering mediums, or the combination of the two with the obstacles as inclusions.

The following is a remark on the geometric and
topological assumptions of the admissible scatterers concerning our
subsequent uniqueness study.

\begin{Remark}
The $C^{0,1}$ and $C^2$ regularity assumptions, respectively, on the refractive index $q$ in $iii)$ and on $\partial\Omega_{ext}$ of the scatterer in $v)$ are only needed for the subsequent uniqueness theorem in determining the location and shape of a scatterer. Though at certain point, such regularity requirement can be weakened, we choose
to work with a consistent assumption to ease the exposition. The topological assumption in $vi)$ is only needed for the
subsequent uniqueness theorem in determining the scattering medium $q$ provided the buried obstacle
have been recovered.
\end{Remark}

Next, we consider the direct scattering problem with a scatterer $\Sigma(\mathbf{D},q,\Omega)\in \mathcal{C}$. Starting from now on, we fix the
wave number to be $\kappa_0>0$. Let $\mathcal{L}_q:=\Delta+\kappa_0^2 q$ denote the Schr\"odinger operator. Using the fact that $(\Delta+\kappa^2)u^i=0$, the forward scattering
problem is reformulated as
\begin{equation}\label{eq:forward problem}
\mathcal{L}_q u^s=f_q(u^i)\quad\mbox{in\ $\mathbf{G}$},\quad \mathcal{B}u^s=g(u^i)\,\,
\mbox{on\ $\partial \mathbf{G}$\,\, and \ $\mathcal{M}u=0$},
\end{equation}
where $f_q(u^i)=\kappa_0^2(1-q)u^i$ and $g(u^i)=-u^i|_{\partial
\mathbf{G}}$. In the sequel, for each
$\rho>0$, we set $\mathbf{G}_\rho:=\mathbf{G}\cap B_\rho(0)$. To
study the direct scattering problem, it is convenient to introduce
the notation
\[
H_{loc}^1(\mathbf{G})=\{u\in \mathscr{D}'(\mathbf{G});
u|_{\mathbf{G}_\rho}\in H^1(\mathbf{G}_\rho)\ \mbox{for each
$\rho>0$ such that $\Sigma\subset B_\rho(0)$}\}.
\]
Then, one can show the well-posedness of the direct scattering
problem in the space $H_{loc}^1(\mathbf{G})$. In fact, the
uniqueness is easily derived by using the Rellich uniqueness theorem (see
Lemma~6.1 in \cite{Isa}). For the existence, by using the
Lax-Phillips method, one can reduce the problem to a corresponding
one in the bounded domain $\Omega_\rho$ (see Chapter 6 in
\cite{Isa}), which has been well understood (see \cite{GilTru} and \cite{Mcl}).
Moreover, we have the well-known elliptic stability estimate, that is,
$\|u\|_{H^1(\mathbf{G}_\rho)}\lesssim
\|f_q(u^i)\|_{L^2(\Omega)}+\|g(u^i)\|_{H^{1/2}(\partial \mathbf{G})}$. However, for our
subsequent uniqueness study in the inverse problem, we need an integral representation
of the solution, which has been given in \cite{KirKre}. To this end, we let
$\Phi(x,y):=e^{i\kappa_0|x-y|}/|x-y|$ be the fundamental solution to the differential operator $(\Delta+\kappa_0^2)$
and introduce the following potential operators:
\begin{equation}
\mbox{SL}\psi(x)=\int_{\partial\mathbf{G}}\Phi(x,y)\psi(y) d S_y,\qquad \mbox{DL}\psi(x)=\int_{\partial\mathbf{G}}
\frac{\partial\Phi(x,y)}{\partial\nu(y)}\psi(y)d S_y
\end{equation}
where $\nu$ is the interior normal to $\mathbf{G}$, and
\begin{equation}
V_q\psi(x)=\kappa_0^2\int_{\Omega}\Phi(x,y)[1-q(y)]\psi(y) d S_y.
\end{equation}
$\mbox{SL}$ and $\mbox{DL}$ are well-known as the single- and double-layer potential operators, while $V_q$
is known as the volume potential operator; we refer \cite{ColKre} and \cite{Mcl} for a detailed study and relevant mapping
properties. For the forward scattering problem, we have the following theorem which is readily modified from Theorem 2.2 in
\cite{KirKre}.
\begin{Theorem}\label{thm:direct}
$u^s\in H_{loc}^1(\mathbf{G})\cap C(\bar{\mathbf{G}})$ is a solution of (\ref{eq:forward problem}) if $u^s|_{\bar{\Omega}}
\in C(\bar{\Omega})$ has the form
\begin{equation}\label{eq:integral 1}
u^s(x)=-V_q u^s(x)+(\mbox{\emph{DL}}+i\kappa_0 \mbox{\emph{SL}})\psi(x)+r(x)\quad x\in \bar{\Omega},
\end{equation}
where $r(x):=-V_q u^i(x)$ and $\psi(x)\in C(\partial\Omega)$ satisfies
\begin{equation}\label{eq:integral 2}
\psi(x)=2\mathscr{T} V_q(x)-2(\mathscr{T} \mbox{\emph{DL}}+i\kappa_0\mathscr{T}\mbox{\emph{SL}})\psi(x)+t(x)\quad x\in \partial\mathbf{G}
\end{equation}
where $t(x):=2\mathscr{T} V_q u^i(x)$ and $\mathscr{T}$ is the one-sided trace operator for $\mathbf{G}$. Moreover, we have
\begin{enumerate}
\item[i)] The system (\ref{eq:integral 1})-(\ref{eq:integral 2}) of integral equations is uniquely solvable in $C(\bar{\Omega})\times C(\partial\Omega)$ for $(r,t)\in C(\bar{\Omega})\times C(\partial\Omega)$ and depends continuously on $r$ and $t$.

\item[ii)] The system (\ref{eq:integral 1})-(\ref{eq:integral 2}) of integral equations is uniquely solvable in $L^2({\Omega})\times C(\partial\Omega)$ for $(r,t)\in L^2({\Omega})\times C(\partial\Omega)$ and depends continuously on $r$ and $t$.
\end{enumerate}
\end{Theorem}

\noindent It is remarked that in our uniqueness study of determining the obstacle, we would essentially make use the continuity of
scattered field in the exterior domain. Next we introduce a more singular point source than $\Phi(x,y)$ which
is given for every fixed $x_0\in \mathbb{R}^3$ by
\begin{equation}
\Psi(y,x_0)=h_1^{(1)}(\kappa_0\rho)P_1(\cos(\psi)),
\end{equation}
where $h_1^{(1)}$ is the spherical Hankel function of the first kind
of order one and  $P_1$ is the Legendre polynomial of order one; and
$(\rho, \phi, \psi)$ is the spherical coordinate of $y-x_0$.
$\Psi(y,x)$ is known as the spherical wave function and we refer to
\cite{ColKre} for related study. It is noted that $\Psi(y,x_0)$ has
quadratic singularity only at the point $y=x_0$ which comes from
that of the spherical Hankel function; that is,
$(y-x_0)^2\Psi(y,x_0)$ is smooth over $\mathbb{R}^3$.

We conclude this section with an approximation property of point sources by linear combination of plane waves.

\begin{Lemma}\label{lem:approximation}
Let $E\subset\mathbb{R}^3$ be a compact set and $x_0\in \mathbb{R}^3\backslash E$ be fixed. Then there exist sequences $v_n(y)$
and $\omega_n(y)$ in the span of plane waves
\[
\mathscr{E}:=\mbox{span}\{e^{i\kappa_0 y \cdot \theta}: \theta\in \mathbb{S}^2\}
\]
such that
\begin{equation}
\|v_n-\Phi(\cdot, x_0)\|_{C^1(E)}\rightarrow 0\quad \mbox{as\ $n\rightarrow \infty$}.
\end{equation}
and
\begin{equation}
\|\omega_n-\Psi(\cdot, x_0)\|_{C^1(E)}\rightarrow 0\quad \mbox{as\ $n\rightarrow \infty$}.
\end{equation}
\end{Lemma}

\begin{proof}
This follows from Lemma~3.2 in \cite{KirKre} by noting that $\Phi(\cdot,x_0)$ and $\Psi(\cdot, x_0)$ are smooth solutions
for the Helmholtz equation in any domain that does not contain $x_0$. See also Lemma~5 in \cite{Pot}.
\end{proof}

\section{Recovery of inhomogeneities together with buried obstacles}\label{sect:uniq with soft obstacle}

In this section, we show the uniqueness in determining a scatterer $\Sigma(\mathbf{D},q,\Omega)\in \mathcal{C}$
by its corresponding scattering amplitude $\mathcal{A}(\hat{x},\theta,\kappa_0)$. The main result is stated as follows.

\begin{Theorem}\label{thm:main}
$\Sigma(\mathbf{D},q,\Omega)\in \mathcal{C}$ is
uniquely determined by knowledge of the far-field pattern
$\mathcal{A}(\hat{x},\theta,\kappa_0)$ for arbitrarily fixed
$\kappa_0>0$ and all $\hat{x},\theta\in \mathbb{S}^{2}$.
\end{Theorem}

The proof of Theorem~\ref{thm:main} is proceeded in three steps, which
we shall outline briefly in the following. In the first step, we recover
the exterior boundary of the scatterer, namely $\partial\Sigma$, disregarding the
interior medium and obstacle. This is based on the use of the singular sources $\Psi(\cdot, x_n)$
with $x_n$ approaching a point $x_0$ which may either lies on the exterior boundary part of the medium
or the exterior boundary part of the obstacle. It is shown the corresponding scattered waves will blow up
in the limiting case. We would like to remark that the use of point source with quadratic singularity is also
considered in \cite{Pot} to determine the support of a scattering medium.
In the second step, we show that one can distinguish the exterior medium boundary from the exterior
obstacle boundary, and thus can determine the inside obstacle. This is based on the use of singular source
$\Phi(\cdot,x_n)$ with $x_n$ approaching an exterior boundary point $x_0$ of the scatterer. It is shown that the scattered
wave will blow up in the limiting case when $x_0$ lies on the boundary of the obstacle, whereas scattered wave remains bounded when
$x_0$ lies on the boundary of the medium. In the final step, we recover the medium along the line of the Sylvester-Uhlmann methodology
(see \cite{SU}). In doing this, we first derive a novel approximation result of Runge type (see Lemma~\ref{lem:3}).
The result is remarkable since it enables us to use CGO solutions in different medium components with different complex phases.
Next, it is natural to construct the almost complex exponential solutions which vanish on the interior boundary of the obstacle.
Since the medium and the buried obstacle have only planar contact, this can be carried out by making reflections
of the CGO solutions with respect to the contact planes. In \cite{Isa3}, similar idea of implementing reflection of solutions
have been used to prove a uniqueness in inverse conductivity problem with local Cauchy data on the boundary. We would like to note
that in \cite{Isa3} the inaccessible part of the boundary is assumed to be on a single plane.

\subsection{Unique determination of $\Sigma$}\label{subsect:3.1}

By contradiction, let
$\tilde{\Sigma}:=\tilde{\Sigma}(\tilde{\mathbf{D}},\tilde{q},\tilde{\Omega})\in
\mathcal{C}$ be a scatterer such that $\tilde{\Sigma}\neq \Sigma$
and
$\mathcal{A}(\hat{x},\theta)=\tilde{\mathcal{A}}(\hat{x},\theta)$
for $\hat{x}, \theta\in \mathbb{S}^2$, where $\mathcal{A}$ and
$\tilde{A}$ are respectively the scattering amplitudes of $\Sigma$
and $\tilde{\Sigma}$ corresponding to the incident plane waves
$\exp\{{i\kappa_0 x\cdot \theta}\}$. Let $\Lambda$ be the (unique)
unbounded connected component of $\mathbb{R}^3\backslash
{(\Sigma\cup \tilde{\Sigma})}$. We denote by
$u(x,\theta)$ and $\tilde{u}(x,\theta)$, respectively, the total
fields corresponding to $\Sigma$ and $\tilde{\Sigma}$. Then, by
the Rellich uniqueness theorem, we know
$u(x,\theta)=\tilde{u}(x,\theta)$ in $\Lambda$ for all $\theta\in
\mathbb{S}^2$. Since $\Sigma\neq \tilde{\Sigma}$ and both are
connected, we easily see that either
$(\mathbb{R}^3\backslash\bar{\Lambda})\backslash{\Sigma}\neq
\emptyset$ or
$(\mathbb{R}^3\backslash\bar{\Lambda})\backslash{\tilde\Sigma}\neq
\emptyset$. Without loss of generality, we assume the former case
and set
$\Sigma^*=(\mathbb{R}^3\backslash\bar{\Lambda})\backslash{\Sigma}$.
It is obvious that $\partial \Sigma^*\subset \partial \Lambda\cup
\partial \Sigma\subset \partial \tilde{\Sigma}\cup \partial\Sigma$
and $\partial\Sigma^*\backslash \partial\Sigma\neq \emptyset$.
According to Definition~\ref{def:scatterer}, $\partial
\tilde{\Sigma}=\partial \tilde{\mathbf{D}}_{ext}\cup \partial
\tilde{\Omega}_{ext}$. Let $x_0\in \partial\Sigma^*\backslash
\partial\Sigma\subset (\partial \tilde{\mathbf{D}}_{ext}\cup \partial
\tilde{\Omega}_{ext})\backslash \Sigma$. We next distinguish
two cases that $x_0\in \partial\tilde{\mathbf{D}}_{ext}\backslash
{\Sigma}$ and $x_0\in \partial \tilde{\Omega}_{ext}\backslash
{\Sigma}$. In the following, we fix $\rho_0>0$ be sufficiently
large such that ${\Sigma\cup\tilde{\Sigma}}\subset
B_{\rho_0}(0)$ and let $\mathbf{G}_{\rho_0}$ and
$\tilde{\mathbf{G}}_{\rho_0}$, respectively, denote
$(\mathbb{R}^3\backslash \bar{\mathbf{D}})\cap B_{\rho_0}(0)$ and
$(\mathbb{R}^3\backslash \bar{\tilde{\mathbf{D}}})\cap B_{\rho_0}(0)$.

\medskip

$Case~1.$~$x_0\in \partial\tilde{\mathbf{D}}_{ext}\backslash
{\Sigma}$. Let $\tau_0>0$ be sufficiently small such that
$B_{\tau_0}(x_0)\subset \Sigma_\infty$ and $B_{\tau_0}(x_0)\cap
\tilde{\Omega}=\emptyset$. Set $S:=\partial \tilde{\mathbf{D}}_{ext}\cap
B_{\tau_0}(x_0)$. Without loss of generality, we assume that
$S\subset \partial\tilde{\mathbf{D}}_{ext}\cap \partial \Lambda$.
Obviously, $B_{\tau_0}(x_0)$ is divided by $S$ into two parts and
we denote by $B_{\tau_0}^+$ the one contained in $\Lambda$. We
now consider the two scattering problems corresponding to
$\Sigma(\mathbf{D},q,\Omega)$ and
$\tilde{\Sigma}:=\tilde{\Sigma}(\tilde{\mathbf{D}},\tilde{q},\tilde{\Omega})$
with the incident fields being the point sources $\Psi(\cdot,x)$ for $x\in B_{\tau_0}^+$.
Let $\omega^s(\cdot,x)$ and $\tilde{\omega}^s(\cdot,x)$ denote,
respectively, the scattered fields. Since the scattered waves coincide in
$\Lambda$ for all plane waves, by using Lemma~\ref{lem:approximation}, it is
straightforward to show that $\omega^s(\cdot, x)=\tilde{\omega}^s(\cdot,x)$ in $\Lambda$
for $x\in B_{\tau_0}^+$. Next, it is observed that $B_{\tau_0}^+\subset
\Sigma_\infty$, and hence $\Psi(\cdot, x)$ with $x\in B_{\tau_0}^+$ is
smooth in $\bar{\Sigma}$. By using the expansions of the spherical
Hankel functions, one can verify directly that
$\|\Psi(\cdot,x)\|_{C^1({\Sigma})}\leq C$ for $x\in
B_{\tau_0}^+$. From the well-posedness of the forward scattering
problem, we see
$\|\omega^s(\cdot,x)\|_{C(\mathbf{G}_{\rho_0})}\leq C$;
see related discussion in Section~\ref{sect:direct scattering}.

Then, we choose $h>0$ such that the sequence
\begin{equation}\label{eq:sequence}
x_n:=x_0+\frac h n \nu(x_0)\qquad n=1,2,\ldots
\end{equation}
is contained in $B_{\tau_0}^+$, where $\nu(x_0)$ is the outward normal to $\partial \tilde{\mathbf{D}}$ at $x_0$.
By our discussion made earlier, $|\omega^s(x_0,x_n)|\leq C$ uniformly for $n\geq 1$.
On the other hand, referring to Lemma~3 in \cite{Pot}, we know $|\Psi(x_0,x_n)|\rightarrow \infty$ as $n\rightarrow\
\infty$, and hence by using the Dirichlet boundary condition of $\tilde{\omega}^s$ on $\partial \tilde{\mathbf{D}}$
\begin{equation}\label{eq:contradiction}
|\omega^s(x_0,x_n)|=|\tilde{\omega}^s(x_0,x_n)|=|-\Psi(x_0,x_n)|\rightarrow\infty\quad \mbox{as\ $n\rightarrow \infty$}.
\end{equation}
This obviously gives a contradiction.

\medskip

$Case~2.$~$x_0\in \partial\tilde{\Omega}_{ext}\backslash
{\Sigma}$. Similar to Case 1, we let $B_{\tau_0}(x_0)$ be a
sufficiently small ball such that
$B_{\tau_0}(x_0)\subset\Sigma_\infty$ and $B_{\tau_0}(x_0)\cap
\bar{\tilde{\mathbf{D}}}=\emptyset$. Moreover, let
$S:=\partial\Omega_{ext}\cap B_{\tau_0}(x_0)$ which is assumed to
lie entirely on $\partial \Lambda$, and let $B_{\tau_0}^+$ denote
the part of $B_{\tau_0}(x_0)$ contained in $\Lambda$. By a same
argument as that for Case~1, we know
$\|\omega^s(\cdot,x)\|_{C({\mathbf{G}}_{\rho_0})}\leq C$ for $x\in
B_{\tau_0}^+$. Clearly, in order to get a contradiction, we only
need to show that $\tilde{\omega}^s(\cdot,x)$ reveals singular
behavior near $x_0$. To this end, let $x_n, n=1,2,\ldots$ be as
defined in (\ref{eq:sequence}). It is first observed that
$|V_{\tilde q}\Psi(x,x_n)|\lesssim 1/|x-x_n|$ (see Lemma~4 in
\cite{Pot}). Hence $\|V_{\tilde
q}\Psi(\cdot,x_n)\|_{L^2(\tilde{\mathbf{G}}_{\rho_0})}\leq C$
uniformly for $n\in \mathbb{N}$. Moreover, noting $x_n$'s are
contained in $B_{\tau_0}^+$ which is away from
$\bar{\tilde{\mathbf{D}}}$, $\|\mathscr{T}V_{\tilde q}
\Psi(\cdot,x_n)\|_{C(\partial\tilde{\mathbf{D}})}\leq C$ uniformly
for $n\in \mathbb{N}$. By Theorem~\ref{thm:direct}, ii),
$\|\tilde{\omega}^s(\cdot,x_n)\|_{L^2(\Omega)}\leq C$ and
$\|\psi(\cdot,x_n)\|_{C(\partial\Omega)}\leq C$, where
$\psi(\cdot,x_n)$ is the density in (\ref{eq:integral 2})
corresponding to the incident waves $\Psi(\cdot,x_n)$. Next, using
the mapping properties that $V_{\tilde q}$ maps $L^2(\Omega)$
continuously into $C({\tilde{\mathbf{G}}}_{\rho_0})$, and
$\mbox{SL}\psi$ and $\mbox{DL}\psi$ map $C(\partial\Omega)$
continuously into $C(\tilde{\mathbf{G}}_{\rho})$ (see
\cite{ColKre}), we know $|V_{\tilde
q}\tilde{\omega}^s(x_0,x_n)|\leq C$ and $|(\mbox{DL}+i\kappa_0
\mbox{SL})\psi(x_0,x_n)|\leq C$ uniformly for $n\in \mathbb{N}$.
On the other hand, referring to Lemma~3 in \cite{Pot}, we know
\begin{equation}\label{eq:split 1}
|V_{\tilde q}\Psi(x_0,x_n)|\rightarrow\infty\quad \mbox{as\ $n\rightarrow\infty$}.
\end{equation}
Hence, by using the relation given in (\ref{eq:integral 1})
\[
|\tilde{\omega}^s(x_0,x_n)|\geq |V_{\tilde q}\Psi(x_0,x_n)|-|V_{\tilde q}\tilde{\omega}(x_0,x_n)|
-|(\mbox{DL}+i\kappa_0 \mbox{SL})\psi(x_0,x_n)|\rightarrow\infty
\]
as\ $n\rightarrow\infty$,
which then yields a similar contradiction to that in
(\ref{eq:contradiction}). \hfill $\Box$

\subsection{Recovery of the obstacle $\mathbf{D}$ }\label{subsect:3.2}

Let $\Sigma(\mathbf{D},q,\Omega)$ and
$\tilde{\Sigma}(\tilde{\mathbf{D}},\tilde{q},\tilde{\Omega})$ be
the two scatterers considered in subsection~\ref{subsect:3.1}, we
next show $\mathbf{D}=\tilde{\mathbf{D}}$. Since
$\Sigma=\tilde{\Sigma}$ and both $\Sigma(\mathbf{D},q,\Omega)$ and
$\tilde{\Sigma}(\tilde{\mathbf{D}},\tilde{q},\tilde{\Omega})$
belongs to class $\mathcal{C}$, one only need to show that
$\partial\mathbf{D}_{ext}=\partial\tilde{\mathbf{D}}_{ext}$ which
then implies $\Gamma_{int}=\tilde{\Gamma}_{int}$ and
$\partial\Omega_{ext}=\partial\tilde{\Omega}_{ext}$. In fact, due
to assumptions iv) and vi) in Definition~\ref{def:scatterer}, it
is easily seen that each planar contact $\Gamma_l$ corresponds
uniquely to a $\Omega_l$. Moreover, noting that $\Gamma_l$ is a
simply connected part of some plane, if $\Sigma=\tilde{\Sigma}$,
$\partial\mathbf{D}_{ext}=\partial\tilde{\mathbf{D}}_{ext}$ and
$\partial\Omega_{ext}=\partial\tilde{\Omega}_{ext}$, one must have
$\mathbf{D}=\tilde{\mathbf{D}}$ and $\Omega=\tilde{\Omega}$. Next,
we assume contrarily that $\partial\mathbf{D}_{ext}\neq
\partial\tilde{\mathbf{D}}_{ext}$. Without loss of generality, let
$\partial\tilde{\mathbf{D}}_{ext}\backslash\partial\mathbf{D}_{ext}
\subset\partial\Sigma\backslash\partial\mathbf{D}_{ext}\subset\partial\Omega_{ext}$
be non-void. Fix $x_0\in
\partial\tilde{\mathbf{D}}_{ext}\backslash\partial\mathbf{D}_{ext}$
and take $B_{\tau_0}(x_0)$ be sufficiently small such that
$B_{\tau_0}(x_0)\cap \bar{\tilde{\mathbf{D}}}=\emptyset$. Let
$x_n$ be as defined in (\ref{eq:sequence}) such that $x_n\in
B_{\tau_0}(x_0)\cap\Sigma_\infty, n=1,2,\ldots$. As before, we
consider the scattering problems corresponding to the point
sources $\Phi(\cdot,x_n)$ and denote by $\omega^s(\cdot,x_n)$ and
$\tilde{\omega}^s(\cdot,x_n)$ the scattered fields corresponding
to $\Sigma(\mathbf{D},q,\Omega)$ and
$\tilde{\Sigma}(\tilde{\mathbf{D}},\tilde{q},\tilde{\Omega})$,
respectively. Obviously, by Lemma~\ref{lem:approximation}, we see
$\omega^s(\cdot,x_n)=\tilde{\omega}^s(\cdot,x_n)$ over
$\Sigma_\infty$.

Since $\|\Phi(\cdot,x_n)\|_{L^2(\Omega)}\leq C$ uniformly for
$n\in \mathbb{N}$, we know
$\|V_{{q}}\Phi(\cdot,x_n)\|_{C(\bar{\mathbf{G}})_{\rho_0}}\leq C$
uniformly for $n\in\mathbb{N}$. By the well-posedness of the
direct scattering problem, $|\omega^s(x_0,x_n)|\leq C$ uniformly
for $n\in \mathbb{N}$. On the other hand, noting $x_0\in
\partial\tilde{\mathbf{D}}$, we have from the homogeneous
Dirichlet boundary condition on $\partial\tilde{\mathbf{D}}$ that
$\tilde{\omega}^s(x_0,x_n)=-\Phi(x_0,x_n)$. Finally, we can get a
contradiction by observing that $|\Phi(x_0,x_n)|\rightarrow\infty$
as $n\rightarrow\infty$. Therefore,
$\mathbf{D}=\tilde{\mathbf{D}}$, which in turn implies
$\Omega=\tilde{\Omega}$. \hfill $\Box$

\subsection{Unique determination of the scattering medium $q$}

In view of the results in subsection~\ref{subsect:3.1} and
\ref{subsect:3.2}, we only need to show that if
$\Sigma(\mathbf{D},q,\Omega)$ and
$\Sigma(\mathbf{D},\tilde{q},\Omega)$ produce the same scattering
amplitude, then $q=\tilde{q}$.

Let $\mathrm{B}:=B_{\rho}(0)$ with suitably selected $\rho>0$ such
that ${\Sigma}\subset \mathrm{B}$, and $\kappa_0^2$ is not a
Dirichlet eigenvalue for $-\Delta$ neither in $\mathrm{B}$.
Moreover, we require that the homogeneous Dirichlet problem for
$\mathcal{L}_{\tilde q}$ has only trivial solution in
$H_0^1(\mathrm{B}\backslash\bar{\mathbf{D}})$. Setting
$w=u-\tilde{u}$, we see
\begin{equation}\label{eq:step I 1}
w=0\qquad \mbox{in\ $\mathrm{B}\backslash {\Sigma}$},
\end{equation}
and hence
\begin{equation}\label{eq:step I 2}
w=\frac{\partial w}{\partial \nu}=0\qquad \mbox{on \  $\partial
\Omega_{ext}:=\partial \Omega\backslash \Gamma_{int}$},
\end{equation}
where $\nu$ is unit outward normal to $\partial \Omega$. Moreover,
by noting that $\mathbf{D}$ is a sound-soft obstacle,
\begin{equation}\label{eq:step I 3}
w=0\qquad \mbox{on\ $\partial \mathbf{D}:=\partial
\mathbf{D}_{ext}\cup\Gamma_{int}$}.
\end{equation}
It is also straightforward to verify that $u\in H^1(\Omega)$
satisfies the following differential equation
\begin{equation}\label{eq:step I 4}
\Delta w+\kappa_0^2 q w=\kappa_0^2 \delta_q \tilde{u},
\end{equation}
where $\delta_q=q-\tilde{q}$. Next, we define
\begin{equation}\label{eq:step I 5}
\mathscr{H}_{q, \Gamma_{int}}:=\{v\in H^{1}(\Omega);
\mathcal{L}_{q} v=0\ \ \mbox{in $\Omega$\ \ and $v=0$\ \ on
$\Gamma_{int}$}\}.
\end{equation}
Multiplying both sides of (\ref{eq:step I 4}) by an arbitrary
$v\in \mathscr{H}_{q, \Gamma_{int}}$ and using Green's formula, we
have
\[
\int_{\Omega} \kappa_0^2 \delta_q \tilde{u} v \ dx\\
= \int_{\Omega} (\mathcal{L}_{q} w)v-(\mathcal{L}_{q} v) w \ dx\\
= \int_{\partial \Omega}\frac{\partial w}{\partial \nu}
v-w\frac{\partial v}{\partial \nu}\ dS_x .
\]
In terms of the relations in (\ref{eq:step I 2})-(\ref{eq:step I
3}), this further yields
\begin{equation}\label{eq:step I 6}
\int_{\Omega} \kappa_0^2 \delta_{q} \tilde{u} v \ dx=0.
\end{equation}
Equivalently, (\ref{eq:step I 6}) is read as
\begin{equation}\label{eq:step I 7}
\sum_{l=0}^{N_0}\int_{\Omega_l}\delta_{q} \tilde{u} v \ dx=
\int_{\Omega} \delta_{q} \tilde{u} v \ dx=0,
\end{equation}
where it is recalled that $\Omega_0$ is separated from the
inhomogeneity while $\Omega_l$ has planar contact with the inhomogeneity
at $\Gamma_l$ for $l=1,2,\ldots, N_0$. We next divide our
argument into three steps.

\medskip

\noindent{\bf Step~I. Denseness argument and two approximation results}

\medskip

Define
\begin{equation}\label{eq:step II 1}
\mathscr{H}_{\tilde q, \mathrm{B}\backslash
\bar{\mathbf{D}}}:=\{\phi\in H^{1}(\mathrm{B}\backslash
\bar{\mathbf{D}}); \mathcal{L}_{\tilde q} \phi=0\ \ \mbox{in
$\mathrm{B}\backslash \bar{\mathbf{D}}$\ \ and $\phi=0$\ \ on
$\partial\mathbf{D}$}\}.
\end{equation}
and
\begin{equation}\label{eq:h q2 gamma}
\mathscr{H}_{\tilde q,\Gamma_{int}}:=\{\phi\in H^1(\Omega);
\mathcal{L}_{\tilde q} \phi=0\ \mbox{in $\Omega$}\ \
\mbox{$\phi=0$\ \ on $\Gamma_{int}$} \}.
\end{equation}

We shall show the following two lemmata at the end of the present
subsection.

\begin{Lemma}\label{lem:1}
The set of total fields $\{u(x; \theta, \tilde{q}); \theta\in
\mathbb{S}^{2}\}$ to (\ref{eq:direct system}) is complete in
$\mathscr{H}_{\tilde q, \mathrm{B}\backslash \bar{\mathbf{D}}}$
with respect to the
$L^2(\mathrm{B}\backslash\bar{\mathbf{D}})$-norm.
\end{Lemma}

\begin{Lemma}\label{lem:2}
Any $\phi\in \mathscr{H}_{\tilde q, \Gamma_{int}}$ can be
$L^2(\Omega)$-approximated by distributions in
$\mathscr{H}_{\tilde q, \mathrm{B}\backslash \bar{\mathbf{D}}}$.
\end{Lemma}

Combining Lemmata~\ref{lem:1} and \ref{lem:2}, we easily see

\begin{Lemma}\label{lem:3}
The set of total fields $\{u(x; \theta, \tilde{q}); \theta\in
\mathbb{S}^{2}\}$ to (\ref{eq:direct system}) is complete in
$\mathscr{H}_{\tilde q, \Gamma_{int}}$ with respect to the
$L^2(\Omega)$-norm.
\end{Lemma}

By Lemma~\ref{lem:3}, we have from (\ref{eq:step I 7}) that
\begin{equation}\label{eq:step II 2}
\int_{\Omega} \delta_{q} \phi v \ dx=0,\qquad \forall \phi\in
\mathscr{H}_{\tilde q, \Gamma_{int}}, \forall v\in \mathscr{H}_{q,
\Gamma_{int}}.
\end{equation}

\medskip

\noindent\textbf{Step~II. Construction of the CGO solutions vanishing on the buried boundary of the obstacle}

\medskip

In this step, we construct special complex geometric optics
solutions to the Schr\"odinger operator $\mathcal{L}_p$ with compactly supported $p\in
L^\infty(\mathbb{R}^3)$.  Let $\xi=(\xi_1,\xi_2,\xi_3)\in
\mathbb{R}^3$. We introduce
\[
e(1)=(\xi_1^2+\xi_2^2)^{-1/2}(\xi_1,\xi_2, 0),\ \ e(3)=(0, 0, 1)
\]
and the unit vector $e(2)$ to form a orthonormal basis $e(1), e(2),
e(3)$ in $\mathbb{R}^3$. The coordinate of $x\in \mathbb{R}^3$ in
this basis is denoted by $(x_{1e}, x_{2e}, x_{3e})_e$. It is
observed that
\[
\xi=(\xi_{1e}, 0, \xi_3)_e, \ \xi_{1e}=(\xi_1^2+\xi_2^2)^{1/2}
\]
and
\[
x\cdot y=\sum_{l=1}^3 x_ly_l=\sum_{l=1}^{3}x_{le}y_{le}.
\]
Define
\begin{equation}\label{eq:zeta}
\begin{split}
\zeta(1)=& (\frac{\xi_{1e}}{2}-\tau\xi_3, i |\xi|(\frac 1
4+\tau^2)^{\frac 1 2}, \frac{\xi_3}{2}+\tau \xi_{1e})_e,\\
\zeta^*(1)=& (\frac{\xi_{1e}}{2}-\tau\xi_3, i |\xi|(\frac 1
4+\tau^2)^{\frac 1 2}, -\frac{\xi_3}{2}-\tau \xi_{1e})_e,\\
\zeta(2)=& (\frac{\xi_{1e}}{2}+\tau\xi_3, -i |\xi|(\frac 1
4+\tau^2)^{\frac 1 2}, \frac{\xi_3}{2}-\tau \xi_{1e})_e,\\
\zeta^*(2)=& (\frac{\xi_{1e}}{2}+\tau\xi_3, -i |\xi|(\frac 1
4+\tau^2)^{\frac 1 2}, -\frac{\xi_3}{2}+\tau \xi_{1e})_e,
\end{split}
\end{equation}
where $\tau$ is a positive real number. By straightforward
calculations, one can verify that
\begin{equation}\label{eq:orthogonal}
\zeta(l)\cdot\zeta(l)=\zeta^*(l)\cdot\zeta^*(l)=0\quad l=1,2.
\end{equation}
From the geometric interpretation of the inner product for
vectors in $\mathbb{R}^3$, we further see that for any unitary matrix $U\in
\mathbb{R}^{3\times 3}$
\begin{equation}\label{eq:orthogonal U}
U\zeta(l)\cdot U\zeta(l)=U\zeta^*(l)\cdot U\zeta^*(l)=0\quad l=1,2.
\end{equation}

Next, we construct special CGO solutions in each sub-domain
$\Omega_l, 1\leq l\leq N_0$, of $\Omega$ for $\mathcal{L}_p$. To
ease our exposition, we fix an arbitrary $\Omega_l$ for the
following construction. We denote $p_l$ the restriction of $p$ on
$\Omega_l$. Then, we extend $p_l\in L^\infty(\Omega_l)$ to
$\mathbb{R}^3$ as follows,
\begin{equation}\label{eq:extension}
\hat{p}_l(x)=\begin{cases}
           p\quad &\mbox{for $x\in \Omega_l$},\\
           \mathscr{R}_l p\quad &\mbox{for $x\in
           \mathscr{R}_l\Omega_l$},\\
           0\quad &\mbox{for $x\in \mathbb{R}^3\backslash
           (\Omega_l\cup \mathscr{R}_l\Omega_l)$},
           \end{cases}
\end{equation}
where and in the following, for
a function $f(x), x\in \mathbb{R}^3$, we denote by $\mathscr{R}_l
f(x)=f(\mathscr{R}_l x)$.
That is, $\hat{p}_l\in L^\infty(\mathbb{R}^3)$ is an odd symmetric function
with respect to $\Pi_l$.

Let $U_l\in \mathbb{R}^{3\times 3}$ be a unitary matrix such that
$U_l^T\Pi_l=\{(x_1,x_2,x_3)\in \mathbb{R}^3; x_3=c_l\}$, where $c_l$
is a constant. Because of the relations in (\ref{eq:orthogonal U}),
it is known that there are CGO solutions of the form (see \cite{SU})
\begin{equation}\label{eq:general cgo}
e^{i U_l\zeta(1)\cdot x}(1+\omega_{1,l}),\ \ e^{iU_l \zeta(2)\cdot
x}(1+\omega_{2,l})
\end{equation}
to the equation $\mathcal{L}_{\hat{p}_l}u=0$ in $\mathbb{R}^3$,
where
\begin{equation}\label{eq:limits}
\|\omega_{1,l}\|_{L^2(\mathrm{B}_0)}+\|\omega_{2,l}\|_{L^2(\mathrm{B}_0)}=0\quad
\mbox{as\ $\tau\rightarrow\infty$},
\end{equation}
with $\mathrm{B}_0\subset \mathbb{R}^3$ a ball containing
$\Omega_l\cup \mathscr{R}_l \Omega_l$.

Set
\begin{equation}\label{eq:cgo}
\begin{split}
\psi_1(x)=& e^{i U_l\zeta(1)\cdot x}(1+\omega_{1,l})-e^{i
U_l\zeta(1)\cdot \mathscr{R}_lx}(1+\mathscr{R}_l\omega_{1,l}),\\
\psi_2(x)=& e^{iU_l \zeta(2)\cdot x}(1+\omega_{2,l})-e^{iU_l
\zeta(2)\cdot \mathscr{R}_lx}(1+\mathscr{R}_l\omega_{2,l}).
\end{split}
\end{equation}
We know that $\psi_1, \psi_2\in H^2(\Omega_l\cup
\mathscr{R}_l\Omega_l)$ solve the differential equation
$\mathcal{L}_{\hat{p}_l} \phi=0$, and
\[
\psi_1=\psi_2=0\quad \mbox{on \ $\Gamma_l$}.
\]
Next, we investigate the product of $\psi_1$ and $\psi_2$ for the
subsequent use. It is first observed that
\begin{equation}
e^{i U_l\zeta(1)\cdot x}(1+\omega_{1,l})= e^{i U_l\zeta(1)\cdot
U_lU_l^Tx}(1+\tilde{\omega}_{1,l})=e^{i \zeta(l)\cdot U_l^T
x}(1+\tilde{\omega}_{1,l}),
\end{equation}
where $\tilde{\omega}_{1,l}(U_l^T x)=\omega_{1,l}(x)$. Setting
$y=U_l^T x$ for $x\in \mathbb{R}^3$, we further have
\begin{equation}
e^{i U_l\zeta(1)\cdot x}(1+\omega_{1,l})=e^{i\zeta(1)\cdot
y}(1+\tilde{\omega}_{1,l}(y)).
\end{equation}
In similar manner, we can treat $e^{i U_l\zeta(1)\cdot
\mathscr{R}_lx}(1+\mathscr{R}_l\omega_{1,l}), e^{iU_l \zeta(2)\cdot
x}(1+\omega_{2,l})$ and $e^{iU_l \zeta(2)\cdot
\mathscr{R}_lx}(1+\mathscr{R}_l\omega_{2,l})$ to get
\begin{equation}
\begin{split}
e^{i U_l\zeta(1)\cdot \mathscr{R}_lx}(1+\mathscr{R}_l\omega_{1,l})=&
e^{i \zeta(1)\cdot \hat{\mathscr{R}}_l
y}(1+\hat{\mathscr{R}}_l\tilde{\omega}_{1,l}(y)),\\
e^{iU_l \zeta(2)\cdot x}(1+\omega_{2,l})=& e^{i \zeta(2)\cdot
y}(1+\tilde{\omega}_{2,l}(y)),\\
e^{i U_l\zeta(2)\cdot \mathscr{R}_lx}(1+\mathscr{R}_l\omega_{1,2})=&
e^{i \zeta(2)\cdot \hat{\mathscr{R}}_l
y}(1+\hat{\mathscr{R}}_l\tilde{\omega}_{2,l}(y)),
\end{split}
\end{equation}
where $\hat{\mathscr{R}}_l$ is the reflection with respect to
$U_l^T \Pi_l=\{(y_1,y_2,y_3)\in \mathbb{R}^3; y_3=c_l\}$. We
remind that obviously $\hat{\mathscr{R}}_l (y_1, y_2, y_3)=(y_1,
y_2, 2c_l-y_3)$. In the following, we denote
$\hat{\mathscr{R}}_ly$ by $y^*$ and $\hat{\mathscr{R}}_l
\tilde{\omega}_{\sigma,l}(y)$ by $\tilde{\omega}^*_{\alpha,l}$,
$\alpha=1,2$. Now,
\begin{equation}\label{eq:product}
\begin{split}
\psi_1(x)\psi_2(x)=& [e^{i\zeta(1)\cdot
y}(1+\tilde{\omega}_{1,l})-e^{i\zeta(1)\cdot
y^*}(1+\tilde{\omega}^*_{1,l})]\\
&\times [e^{i\zeta(2)\cdot
y}(1+\tilde{\omega}_{1,2})-e^{i\zeta(2)\cdot
y^*}(1+\tilde{\omega}^*_{1,2})]\\
=& e^{i(\zeta(1)+\zeta(2))\cdot
y}(1+\tilde{\omega}_{1,l})(1+\tilde{\omega}_{2,l})\\
&- e^{i(\zeta(1)+\zeta^*(2))\cdot y}
e^{i\zeta(2)\cdot(0,0,2c_l)}(1+\tilde{\omega}_{1,l})(1+\tilde{\omega}_{2,l}^*)\\
&- e^{i(\zeta(1)^*+\zeta(2))\cdot
y}e^{i\zeta(1)\cdot(0,0,2c_l)}(1+\tilde{\omega}_{1,l}^*)(1+\tilde{\omega}_{2,l})\\
&+ e^{i(\zeta(1)+\zeta(2))\cdot
y^*}(1+\tilde{\omega}_{1,l}^*)(1+\tilde{\omega}_{2,l}^*)\\
=& e^{i\xi\cdot y}(1+\tilde{\omega}_{1,l})(1+\tilde{\omega}_{2,l})\\
&- e^{i(\xi_{1e}, 0, 2\tau\xi_{1e})\cdot y} e^{i\zeta(2)\cdot
(0,0,2c_l)}(1+\tilde{\omega}_{1,l})(1+\tilde{\omega}^*_{2,l})\\
&- e^{i(\xi_{1e}, 0,-2\tau\xi_{1e})\cdot
y}e^{i\zeta(1)\cdot(0,0,2c_l)}(1+\tilde{\omega}^*_{1,l})(1+\tilde{\omega}_{2,l})\\
&+e^{i\xi\cdot
y^*}(1+\tilde{\omega}_{1,l}^*)(1+\tilde{\omega}_{2,l}^*),
\end{split}
\end{equation}
where $y=U^T_lx$. Here we note that in (\ref{eq:product})
\begin{equation}\label{eq:real numbers}
\begin{split}
\zeta(1)\cdot (0,0,2c_l)=& (\xi_3+2\tau\xi_{1e})c_l,\\
\zeta(2)\cdot (0,0,2c_l)=& (-\xi_3+2\tau\xi_{1e})c_l,
\end{split}
\end{equation}
both are real numbers.

\medskip

\noindent\textbf{Step~III. Concluding the
proof}\label{subsect:3.3}

\medskip

With the above preparations, we can conclude the proof
as follows.  First, we fix an $\eta=(\eta_1,\eta_2,\eta_3)\in
\mathbb{R}^3$. Next, as in Step~II, we construct CGO solutions for
the operators $\mathcal{L}_{\tilde q}$ and $\mathcal{L}_{q}$,
respectively as $\psi_1$ and $\psi_2$ in (\ref{eq:cgo}), in each
subdomain $\Omega_l$ with $\xi^l:=U^T_l\eta$ replacing the $\xi$ in
(\ref{eq:zeta}) in each
$\Omega_l$, where $U_l^T\in \mathbb{R}^{3\times
3}$ are unitary matrices such that there are constants $c_l$ and
$U_l^T {\Pi}_l=\{(y_1,y_2,y_3)\in \mathbb{R}^3; y_3=c_l\}$
for $l=1,2,\ldots,N_0$. Whereas for CGO solutions
in $\Omega_0$, we take $\xi^0:=\eta$
in (\ref{eq:zeta}) for defining the complex phases and let them be given as those
in (\ref{eq:general cgo}) without the rotation matrix $U_l$.
That is, we need not the rotation of the subdomain $\Omega_0$, nor the reflection
of CGO solutions in $\Omega_0$.
Noting that the sub-domains $\Omega_l$'s of
$\Omega$ are disjoint from each other, these CGO solutions constructed in each subdomain are
patched together to yield, respectively solutions $\phi\in
\mathscr{H}_{\tilde q,\Gamma_{int}}$ and $v\in \mathscr{H}_{q,\Gamma_{int}}$.
Then, in view of (\ref{eq:step II 2}) and (\ref{eq:product}), we
have
\begin{equation}\label{eq:conluding one}
\begin{split}
0=&\int_{\Omega}\delta_{q}\phi v\
dx=\int_{\Omega_0}\delta_q^0 \phi v\ dx+\sum_{l=1}^{N_0}\int_{\Omega_l}\delta_{q}^l \phi v\ dx\\
=& \int_{\Omega_0}\delta_q^0 e^{i\eta\cdot
x}(1+\omega_{1,0})(1+\omega_{2,0})\ dx+
\sum_{l=1}^{N_0}\int_{\Omega_l}\delta_{q}^l[e^{i\xi^l\cdot
y}(1+\tilde{\omega}_{1,l})(1+\tilde{\omega}_{2,l})\\
&- e^{i(\xi_{1e}^l, 0, 2\tau\xi_{1e}^l)\cdot y} e^{i\zeta(2)\cdot
(0,0,2c_l)}(1+\tilde{\omega}_{1,l})(1+\tilde{\omega}^*_{2,l})\\
&- e^{i(\xi_{1e}^l, 0,-2\tau\xi_{1e}^l)\cdot
y}e^{i\zeta(1)\cdot(0,0,2c_l)}(1+\tilde{\omega}^*_{1,l})(1+\tilde{\omega}_{2,l})\\
&+e^{i\xi^l\cdot
y^*}(1+\tilde{\omega}_{1,l}^*)(1+\tilde{\omega}_{2,l}^*)]\ dx,
\end{split}
\end{equation}
where $\delta_{q}^l$ is the restriction of $\delta_{q}$ on
$\Omega_l$. Clearly, the moduli of all exponents are bounded by
$1$ by noting (\ref{eq:real numbers}). Now, we let
$\tau\rightarrow\infty$ in (\ref{eq:conluding one}). Due to
(\ref{eq:limits}), the limits of all terms containing
$\omega_{\sigma,0}$, $\tilde{\omega}_{\sigma,l}$ and
$\tilde{\omega}_{\sigma,l}^*$, $\sigma=1,2$, are zero. By the
Riemann-Lebsegue Lemma,
\begin{align*}
&\lim_{\tau\rightarrow\infty}\int_{\Omega_l} \delta_{q}^l
e^{i(\xi_{1e}^l, 0, 2\tau\xi_{1e}^l)\cdot y} e^{i\zeta(2)\cdot
(0,0,2c_l)}\ dx\\
=&\lim_{\tau\rightarrow\infty}\int_{\Omega_l}\delta_{q}^l
e^{i(\xi_{1e}^ly_{1e}+\xi_3^l c_l+2\tau\xi_{1e}^l(y_3-c_l))}\ d x=0,\\
&\lim_{\tau\rightarrow\infty}\int_{\Omega_l}\delta_{q}^l
e^{i(\xi_{1e}^l,
0,-2\tau\xi_{1e}^l)\cdot y}e^{i\zeta(1)\cdot(0,0,2c_l)}\ dx\\
=&\lim_{\tau\rightarrow\infty}\int_{\Omega_l}\delta_{q}^l
e^{i(\xi_{1e}^ly_{1e}+\xi_3^l c_l-2\tau\xi_{1e}^l(y_3-c_l))}\ d
x=0,
\end{align*}
provided $\xi_{1e}^l\neq 0$. Define
\[
\mathscr{A}:=\{\eta=(\eta_1,\eta_2,\eta_3)\in \mathbb{R}^3;\ \
\xi_{1e}^l\neq 0\ \mbox{with}\ \xi^l=U^T_l \eta\ \mbox{for}\
l=1,2,\ldots,N_0\}.
\]
Obviously, $\mathscr{A}$ is an open set in $\mathbb{R}^3$. Now,
summarizing the above discussion by letting $\eta\in \mathscr{A}$,
we have obtained from (\ref{eq:conluding one}) that
\begin{equation}\label{eq:concluding two}
\begin{split}
0=&\sum_{l=1}^{N_0}\int_{\Omega_l} \tilde{q}_l (e^{i\xi^l\cdot
y}+e^{i\xi^l\cdot y^*})\ dx+\int_{\Omega_0}\delta_q^0 e^{i\eta\cdot x}\ dx\\
=&\sum_{l=1}^{N_0}\int_{\Omega_l} \tilde{q}_l(e^{i U_l^T\eta\cdot
U_l^T x}+e^{i U_l^T\eta\cdot U_l^T \mathscr{R}_lx})\
dx+\int_{\Omega_0}\delta_q^0 e^{i\eta\cdot x}\ dx\\
=& \sum_{l=1}^{N_0}\int_{\Omega_l}\tilde{q}_l(e^{i\eta\cdot
x}+e^{i\eta\cdot \mathscr{R}_l x})\ dx+\int_{\Omega_0}\delta_q^0
e^{i\eta\cdot x}\ dx.
\end{split}
\end{equation}
As in Step~II, we extend $\delta_{q}$ to $\mathbb{R}^3$ by
patching together those $\delta_{q}^l$'s in $\Omega_l\cup
\mathscr{R}_l \Omega_l$ which is obtained by even extension of
$\delta_{q}^l$ in $\Omega_l$ with respect to ${\Pi}_l$ (cf.
(\ref{eq:extension})), and letting it be zero in
$(\mathbb{R}^3\backslash \Omega_0)\backslash \bigcup_{l=1}^{N_0}
(\Omega_l\cup\mathscr{R}_l\Omega_l)$. This is possible by our
assumption vi) in Definition~\ref{def:scatterer} that
$(\Omega_l\cup \mathscr{R}_l\Omega_l)\cap (\Omega_{l'}\cup
\mathscr{R}_{l'}\Omega_{l'})=\emptyset$ and $(\Omega_l\cup
\mathscr{R}_l\Omega_l)\cap \Omega_0=\emptyset$ for $1\leq l,
l'\leq N_0$ and $l\neq l'$. Hence, we further have from
(\ref{eq:concluding two}) that
\begin{equation}\label{eq:concluding three}
\int_{\mathbb{R}^3} \delta_{q}(x) e^{i\eta\cdot x} \ dx=0,
\end{equation}
for all $\eta\in \mathscr{A}$. Since $\delta_{q}(x)$ is compactly
supported, the LHS of (\ref{eq:concluding three}) is analytic with
respect to $\eta$. So, we see that (\ref{eq:concluding three})
holds for all $\eta\in \mathbb{R}^3$. Now, $\delta_{q}=0$ by the
uniqueness of inverse Fourier transform. The proof is
completed.\hfill $\Box$

\begin{proof}[Proof of Lemma~\ref{lem:1}]

By contradiction, we assume there exists $\bar{f}\in
L^2(\mathrm{B}\backslash\bar{\mathbf{D}})$ such that
\begin{equation}\label{eq:eq 0}
\int_{\mathrm{B}\backslash\bar{\mathbf{D}}} f(x) u(x;\theta,\tilde{q})\ dx=0
\end{equation}
for all total fields $\tilde{u}(x):=u(x;\theta,\tilde{q})$ to (\ref{eq:direct system}) with
$\theta\in \mathbb{S}^2$; whereas
\begin{equation}\label{eq:neq 0}
\int_{\mathrm{B}\backslash\bar{\mathbf{D}}} f \phi \ dx\neq 0
\end{equation}
for some $\phi \in \mathscr{H}_{\tilde q,\mathrm{B}\backslash
\bar{\mathbf{D}}}$. We extend $f$ to be zero in
$\mathbb{R}^3\backslash\mathrm{B}$. Let $u^*\in
H_{loc}^1(\mathbf{G})$ be the unique solution to
\begin{equation}\label{eq:auxilary}
\begin{cases}
&\mathcal{L}_{\tilde q} u^*(x)=f(x)\qquad x\in \mathbf{G},\\
& \quad u^*=0\hspace*{1.8cm}\mbox{on $\partial \mathbf{D}$},
\end{cases}
\end{equation}
and $\mathcal{M}u^*=0$, namely $u^*$ satisfies the radiation condition. In view of
(\ref{eq:eq 0}) and (\ref{eq:auxilary}), we see
\begin{equation}\label{eq:int1}
\int_{\mathrm{B}\backslash\bar{\mathbf{D}}} \tilde u \mathcal{L}_{\tilde q}
u^*=0.
\end{equation}
By further noting $\mathcal{L}_{\tilde q} u=0$, we have from
(\ref{eq:int1}) with the help of Green's formula that
\begin{equation}\label{eq:int2}
\begin{split}
0=& \int_{\mathrm{B}\backslash \bar{\mathbf{D}}}
(\mathcal{L}_{\tilde q}u^*)\tilde u-u^*(\mathcal{L}_{\tilde q}\tilde u)\ dx\\
=& \int_{\partial \mathrm{B}} \frac{\partial u^*}{\partial
\nu}\tilde u-u^*\frac{\partial \tilde u}{\partial \nu} \ dS_x+\int_{\partial
\mathbf{D}} \frac{\partial u^*}{\partial \nu}\tilde u-u^*\frac{\partial
\tilde u}{\partial \nu}\ dx\\
=& \int_{\partial \mathrm{B}} \frac{\partial u^*}{\partial
\nu}\tilde u-u^*\frac{\partial \tilde u}{\partial \nu} \ dS_x,
\end{split}
\end{equation}
where $\nu$ is the exterior normal to corresponding domains and in
the last equality we have made use of boundary conditions $\tilde u=u^*=0$
on $\partial \mathbf{D}$. Next, using the fact $\tilde u(x,
\theta)=\tilde{u}^s(x,\theta)+u^i(x,\theta)=\tilde{u}^s+e^{i\kappa x\cdot \theta}$,
we have from (\ref{eq:int2})
\begin{equation}\label{eq:int3}
\int_{\partial \mathrm{B}} \frac{\partial u^*}{\partial
\nu}u^i-u^*\frac{\partial u^i}{\partial \nu} \ dS_x=-\int_{\partial
\mathrm{B}} \frac{\partial u^*}{\partial \nu}\tilde{u}^s-u^*\frac{\partial
\tilde{u}^s}{\partial \nu} \ dS_x.
\end{equation}
Since $(\Delta+\kappa_0^2)u^*=(\Delta+\kappa_0^2)\tilde{u}^s=0$ in
$\mathbb{R}^3\backslash \mathrm{B}$ and both $u^*$ and $u^s$
satisfies the radiation condition, we see that the RHS of
(\ref{eq:int3}) vanishes identically, and hence
\begin{equation}\label{eq:int4}
\int_{\partial \mathrm{B}} \frac{\partial u^*}{\partial
\nu}u^i-u^*\frac{\partial u^i}{\partial \nu} \ dS_x=0.
\end{equation}
Then we define $\omega^*$ to be the unique solution to
$(\Delta+\kappa^2_0)\omega^*=0$ in $\mathrm{B}$ with Dirichlet
boundary data $\omega^*=u^*$ on $\partial \mathrm{B}$. It is remarked that the unique
existence is guaranteed by our earlier assumption that $\kappa_0^2$ is not a Dirichlet eigenvalue for
$-\Delta$ in $\mathrm{B}$.  Noting that $(\Delta
+\kappa_0^2)u^i=0$ in $\mathrm{B}$, we have from Green's formula
\begin{equation}
0=\int_{\partial\mathrm{B}}\frac{\partial \omega^*}{\partial
\nu}u^i-\omega^*\frac{\partial u^i}{\partial \nu}\
dS_x=\int_{\partial\mathrm{B}}\frac{\partial \omega^*}{\partial
\nu}u^i-u^*\frac{\partial u^i}{\partial \nu}\ dS_x,
\end{equation}
which together with (\ref{eq:int4}) further yields
\begin{equation}\label{eq:p1}
\int_{\partial \mathrm{B}}(\frac{\partial \omega^*}{\partial
\nu}-\frac{\partial u^*}{\partial \nu}) e^{i\kappa_0x\cdot \theta}\
dx=0\quad \mbox{for all \ $\theta\in \mathbb{S}^2$}.
\end{equation}
Since $\kappa_0^2$ is not a Dirichlet eigenvalue for $-\Delta$ in
$\mathrm{B}$, $\{e^{i\kappa_0x\cdot\theta}|_{\partial\mathrm{B}}; \theta\in \mathbb{S}^2\}$ is dense
in $L^2(\partial\mathrm{B})$ (cf. \cite{Isa2}). Hence, one can conclude from (\ref{eq:p1}) that $\frac{\partial
\omega^*}{\partial \nu}=\frac{\partial u^*}{\partial \nu}$ on
$\partial \mathrm{B}$.  Now, if we set $\Psi$ to be $u^*$ in
$\mathbb{R}^3\backslash \mathrm{B}$ and $\omega^*$ in $\mathrm{B}$, then it is an
entire solution to $\Delta+\kappa_0^2$ and satisfies the radiation
condition as $u^*$ does. Clearly, $\Psi$ must be identically zero.
In doing this, we have shown that $u^*=0$ in $\mathbb{R}^3\backslash
\mathrm{B}$. Finally, again by using Green's formula, we have
\begin{align*}
\int_{\mathrm{B}\backslash \bar{\mathbf{D}}} f \phi\ dx=&
\int_{\mathrm{B}\backslash \bar{\mathbf{D}}}
(\mathcal{L}_{\tilde q}u^*)\phi-u^*(\mathcal{L}_{\tilde q}\phi)\ dx\\
=& \int_{\partial \mathrm{B}\cup \partial \mathbf{D}}\frac{\partial
u^*}{\partial \nu}\phi-u^*\frac{\partial \phi}{\partial \nu}
 =0,
\end{align*}
where in the last equality we have made use of homogeneous boundary
conditions $u^*=\frac{\partial u^*}{\partial \nu}=0$ on $\partial
\mathrm{B}$ and $u^*=\phi=0$ on $\partial \mathbf{D}$. This obviously
contradicts to (\ref{eq:neq 0}), thus completing the proof.
\end{proof}

\begin{proof}[Proof of Lemma~\ref{lem:2}]

We assume contrarily that there exits $\bar{f}\in
L^2(\mathrm{B}\backslash \bar{\mathbf{D}})$ supported in $\Omega$ such that
\begin{equation}\label{eq:int21}
\int_{\Omega} f u\ dx =0\quad \forall u\in
\mathscr{H}_{\tilde q,\mathrm{B}\backslash \bar{\mathbf{D}}},
\end{equation}
but
\begin{equation}\label{eq:int22}
\int_{\Omega} f v\ dx\neq 0\quad \mbox{for some\ $v\in
\mathscr{H}_{\tilde q,\Gamma_{int}}$}.
\end{equation}
Let $u^*\in
H_0^1(\mathrm{B}\backslash\bar{\mathbf{D}})$ be the unique solution
to $\mathcal{L}_{\tilde q}u^*=f$. Here it is noted that the unique existence is guaranteed by our earlier requirement
that $\mathrm{B}$ is chosen such that the homogeneous Dirichlet problem for the partial differential operator
$\mathcal{L}_{\tilde q}$ has only trivial solution in $\mathrm{B}\backslash\bar{\mathbf{D}}$. Then,
in view of (\ref{eq:int21}) and with the help of Green's formula, we
have by straightforward calculations that
\begin{equation}\label{eq:int23}
\begin{split}
0=&\int_\Omega f u \ dx=\int_{\Omega}(\mathcal{L}_{\tilde q}
u^*)u-u^*(\mathcal{L}_{\tilde q}u)\ dx\\
=&\int_{\partial \Omega}\frac{\partial u^*}{\partial
\nu}u-u^*\frac{\partial u}{\partial \nu}\ d S_x
=\int_{(\partial \Omega\backslash \Gamma_{int})\cup \partial\mathbf{D}_{ext}}
\frac{\partial u^*}{\partial \nu}u-u^*\frac{\partial u}{\partial
\nu}\ d S_x\\
=&\int_{\partial\Sigma}\frac{\partial u^*}{\partial \nu}u-u^*\frac{\partial u}{\partial
\nu}\ d S_x=\int_{\partial\mathrm{B}} \frac{\partial u^*}{\partial
\nu}u-u^*\frac{\partial u}{\partial \nu}\ d S_x\\
=& \int_{\partial \mathrm{B}}\frac{\partial u^*}{\partial \nu} u\
dx,
\end{split}
\end{equation}
where $\nu$ is the exterior normal to corresponding domains. In the
above deduction, we have made use of the boundary conditions $u=0$
on $\partial \mathbf{D}=\Gamma_{int}\cup \partial\mathbf{D}_{ext}$ and $u^*=0$ on
$\partial\mathbf{D}\cup \partial \mathrm{B}$. It is clear that $u$
can be arbitrary smooth function on $\partial \mathrm{B}$. So, we
have from (\ref{eq:int23}) that $\frac{\partial u^*}{\partial
\nu}=0$ on $\partial \mathrm{B}$. Hence, by the unique continuation
principle, we know $u^*=0$ in
$(\mathrm{B}\backslash\bar{\Sigma}$. Now,
again by using Green's formula, we have
\begin{align*}
\int_{\Omega}f v\ dx=&\int_{\Omega} (\mathcal{L}_{\tilde q} u^*) v-u^*
(\mathcal{L}_{\tilde q}v)\ dx\\
=&\int_{\partial \Omega}\frac{\partial u^*}{\partial \nu}
v-u^*\frac{\partial v}{\partial \nu}\ dS_x\\
=& \int_{\partial \Omega\backslash \Gamma_{int}} \frac{\partial
u^*}{\partial \nu} v-u^*\frac{\partial v}{\partial \nu}\
dS_x+\int_{\Gamma_{int}} \frac{\partial u^*}{\partial \nu}
v-u^*\frac{\partial v}{\partial \nu}\ dS_x=0,
\end{align*}
where we have made use of the homogeneous boundary conditions
$u^*=\frac{\partial u^*}{\partial \nu}=0$ on $\partial
\Omega\backslash \Gamma_{int}$ and $u^*=v=0$ on $\Gamma_{int}$. This
obviously contradicts to (\ref{eq:int22}), which completes the
proof.

\end{proof}

\section{Recovery of scattering medium with known included
obstacle}\label{sect:concluding remark}

As can be seen from the argument in subsection~\ref{subsect:3.2},
in order to recover the buried obstacle, one has to assume that
the obstacle is partly exposed to the exterior of the medium. In
this final section, we would like to remark an interesting case
that one can recover the surrounding medium even if the obstacle
is buried completely but known \emph{a priori}. We would only give
a simple example though one can appeal for a more general study.

Let $\mathbf{D}$ be a bounded polyhedron in $\mathbb{R}^3$ and
$\mathbf{G}=\mathbb{R}^3\backslash\bar{\mathbf{D}}$. We denote by
$F_l, l=1,2,\ldots, m$ the faces of $\mathbf{D}$. For each $F_l$,
we let $\Omega_l\subset\mathbf{G}$ be a bounded Lipschitz domain
such that $\partial\Omega_l\cap\partial\mathbf{D}=F_l$,
$l=1,2,\ldots, m$. We further assume that all $\Omega_l$'s are
simply connected and satisfy a topological requirement as that
given in vi), Definition~\ref{def:scatterer}. Let $q\in
L^\infty(\mathbf{G})$ such that
$supp(1-q)=\cup_{l=1}^m\bar{\Omega}_l$. Clearly, the obstacle
$\mathbf{D}$ is now completely included in the scattering medium.
For such a scatterer, we would like to remark that by using
Lax-Phillips method, one can still show the unique existence of a
solution $u^s\in H_{loc}^1(\mathbf{G})$ to the forward scattering
problem (\ref{eq:forward problem}). It is also readily seen that
the approximation result in Lemma~\ref{lem:3} still holds. Hence,
all our arguments in subsection~\ref{subsect:3.3} remain valid to
show the unique determination of the scattering medium $q$
provided the obstacle $\mathbf{D}$ is known in advance.

\section*{Acknowledgement}

The author would like to thank Prof. Gunther Uhlmann for proposing
the research project and a lot of stimulating discussion.

\end{document}